\documentclass{amsart}
\usepackage{amsfonts}

\newtheorem{theorem}{Theorem}
\theoremstyle{plain}

\numberwithin{equation}{section}
\input{tcilatex}

\begin{document}
\title{On the General Form of Quantum Stochastic Evolution Equation.}
\author{V. P. Belavkin.}
\address{Mathematics Department, University of Nottingham,\\
NG7 2RD, UK.}
\email{vpb@maths.nott.ac.uk}
\date{July 20, 1995}
\subjclass{Quantum Stochastics}
\thanks{Published in: \textit{Stochastic Analysis and Applications} 91--106,
World Scientific, Singapore, 1996.}
\maketitle

\begin{abstract}
A characterisation of the quantum stochastic bounded generators of
irreversible quantum state evolutions is given. This suggests the general
form of quantum stochastic evolution equation with respect to the Poisson
(jumps), Wiener (diffusion) or general Quantum Noise. The corresponding
irreversible Heisenberg evolution in terms of stochastic completely positive
(CP) cocycles is also characterized and the general form of the stochastic
completely dissipative (CD) operator equation is discovered.
\end{abstract}

\section{Quantum Stochastic Filtering Equations}

The quantum filtering theory, which was outlined in \cite{1, 2} and
developed then since \cite{3}, provides the derivations for new types of
irreversible stochastic equations for quantum states, giving the dynamical
solution for the well-known quantum measurement problem. Some particular
types of such equations have been considered recently in the
phenomenological theories of quantum permanent reduction \cite{4,5},
continuous measurement collapse \cite{6,7}, spontaneous jumps \cite{8,9},
diffusions and localizations \cite{10,11}. The main feature of such dynamics
is that the reduced irreversible evolution can be described in terms of a
linear dissipative stochastic wave equation, the solution to which is
normalized only in the mean square sense.

The simplest dynamics of this kind is described by the continuous filtering
wave propagators $V_{t}\left( \omega \right) $, defined on the space $\Omega 
$ of all Brownian trajectories as an adapted operator-valued stochastic
process in the system Hilbert space $\mathcal{H}$, satisfying the stochastic
diffusion equation 
\begin{equation}
\mathrm{d}V_{t}+KV_{t}\mathrm{d}t=LV_{t}\mathrm{dQ},\quad V_{0}=I
\label{1.3}
\end{equation}%
in the It\^{o} sense, which was derived from a unitary evolution in \cite{13}%
. Here $\mathrm{Q}\left( t,\omega \right) $ is the standard Wiener process,
which is described by the independent increments $\mathrm{dQ}\left( t\right)
=\mathrm{Q}\left( t+\mathrm{d}t\right) -\mathrm{Q}\left( t\right) $, having
the zero mean values $\langle \mathrm{dQ}\rangle =0$ and the multiplication
property $(\mathrm{dQ})^{2}=\mathrm{d}t$, $K$ is an accretive operator, $%
K+K^{\dagger }\geq L^{\dagger }L$, and $L$ is a linear operator $\mathcal{D}%
\rightarrow \mathcal{H}$. Using the It\^{o} formula 
\begin{equation}
\mathrm{d}\left( V_{t}^{\dagger }V_{t}\right) =\mathrm{d}V_{t}^{\dagger
}V_{t}+V_{t}^{\dagger }\mathrm{d}V_{t}+\mathrm{d}V_{t}^{\dagger }\mathrm{d}%
V_{t},  \label{1.2}
\end{equation}%
and averaging $\left\langle \cdot \right\rangle $ over the trajectories of $%
\mathrm{Q}$, one obtains $\mathrm{d}\langle V_{t}^{\dagger }V_{t}\rangle
\leq 0$ as a consequence of $L^{\dagger }L\leq K+K^{\dagger }$. Note that
the process $V_{t}$ is necessarily unitary if the filtering condition $%
K^{\dagger }+K=L^{\dagger }L$ holds, and if $L^{\dagger }=-L$ in the bounded
case.

Another type of the filtering wave propagator $V_t\left( \omega \right)
:\psi _0\in \mathcal{H}\mapsto \psi _t\left( \omega \right) $ in $\mathcal{H}
$ is given by the stochastic jump equation 
\begin{equation}
\mathrm{d}V_t+KV_t\mathrm{d}t=LV_t\mathrm{dP},\quad V_0=I,  \label{1.1}
\end{equation}
derived in \cite{12} by the conditioning with respect to the spontaneous
stationary reductions at the random time instants $\omega =\left\{
t_1,t_2,...\right\} $. Here $L=J-I$ is the jump operator, corresponding to
the stationary discontinuous evolutions $J:\psi _t\mapsto \psi _{t+}$ at $%
t\in \omega $, and $\mathrm{P}\left( t,\omega \right) $ is the standard
Poisson process, counting the number $\left| \omega \cap [0,t)\right| $
compensated by its mean value $t$. It is described as the process with
independent increments $\mathrm{dP}\left( t\right) =\mathrm{P}\left( t+%
\mathrm{d}t\right) -\mathrm{P}\left( t\right) $, having the values $\left\{
0,1\right\} $ at $\mathrm{d}t\rightarrow 0$, with zero mean $\langle \mathrm{%
dP}\rangle =0$, and the multiplication property $\left( \mathrm{dP}\right)
^2=\mathrm{dP}+\mathrm{d}t$. Using the It\^{o} formula (\ref{1.2}) with $%
\mathrm{d}V_t^{\dagger }\mathrm{d}V_t=V_t^{\dagger }L^{\dagger }LV_t(\mathrm{%
dP}\mathbf{+}\mathrm{d}t\mathrm{)}$, one can obtain 
\begin{equation*}
\mathrm{d}\left( V_t^{\dagger }V_t\right) =V_t^{\dagger }\left( L^{\dagger
}L-K-K^{\dagger }\right) V_t\mathrm{d}t+V_t^{\dagger }\left( L^{\dagger
}+L+L^{\dagger }L\right) V_t\mathrm{dP}.
\end{equation*}
Averaging $\left\langle \cdot \right\rangle $ over the trajectories of $%
\mathrm{P}$, one can easily find that $\mathrm{d}\langle V_t^{\dagger
}V_t\rangle \leq 0$ under the sub-filtering condition $L^{\dagger }L\leq
K+K^{\dagger }$. Such evolution is unitary if $L^{\dagger }L=K+K^{\dagger }$
and if the jumps are isometric, $J^{\dagger }J=I$.

This proves in both cases that the stochastic wave function $\psi _t\left(
\omega \right) =V_t\left( \omega \right) \psi _0$ is not normalized for each 
$\omega $, but it is normalized in the mean square sense to the probability $%
\langle ||\psi _t||^2\rangle \leq ||\psi _0||^2=1$ for the quantum system
not to be demolished during its observation up to the time $t$. If $%
\left\langle ||\psi _t||^2\right\rangle =1$, then the positive stochastic
function $||\psi _t\left( \omega \right) ||^2$ is the probability density of
a diffusive $\widehat{\mathrm{Q}}$ or counting $\widehat{\mathrm{P}}$ output
process up to the given $t$ with respect to the standard Wiener $\mathrm{Q}$
or Poisson $\mathrm{P}$ input processes.

Using the It\^{o} formula for $\phi _t\left( B\right) =V_t^{\dagger }BV_t$,
one can obtain the stochastic equations 
\begin{equation}
\mathrm{d}\phi _t\left( B\right) +\phi _t\left( K^{\dagger }B+BK-L^{\dagger
}BL\right) \mathrm{d}t=\phi _t\left( L^{\dagger }B+BL\right) \mathrm{dQ},
\label{1.5}
\end{equation}
\begin{equation}
\mathrm{d}\phi _t\left( B\right) +\phi _t\left( K^{\dagger }B+BK-L^{\dagger
}BL\right) \mathrm{d}t=\phi _t\left( J^{\dagger }BJ-B\right) \mathrm{dP},
\label{1.4}
\end{equation}
describing the stochastic evolution $Y_t=\phi _t\left( B\right) $ of a
bounded system operator $B\in \mathcal{L}\left( \mathcal{H}\right) $ as $%
Y_t\left( \omega \right) =V_t\left( \omega \right) ^{\dagger }BV_t\left(
\omega \right) $. The maps $\phi _t:B\mapsto Y_t$ are Hermitian in the sense
that $Y_t^{\dagger }=Y_t$ if $B^{\dagger }=B$, but in contrast to the usual
Hamiltonian dynamics, are not multiplicative in general, $\phi _t\left(
B^{\dagger }C\right) \neq \phi _t\left( B\right) ^{\dagger }\phi _t\left(
C\right) $, even if they are not averaged with respect to $\omega $.
Moreover, they are usually not normalized, $M_t\left( \omega \right) :=\phi
_t\left( \omega ,I\right) \neq I$, although the stochastic positive
operators $M_t=V_t^{\dagger }V_t$ under the filtering condition are usually
normalized in the mean, $\langle M_t\rangle =I$, and satisfy the martingale
property $\epsilon _t\left[ M_s\right] =M_t$ for all $s>t$, where $\epsilon
_t$ is the conditional expectation with respect to the history of the
processes $\mathrm{P}$ or $\mathrm{Q}$ up to time $t$.

Although the filtering equations (\ref{1.1}), (\ref{1.3}) look very
different, they can be unified in the form of quantum stochastic equation 
\begin{equation}
\mathrm{d}V_{t}+KV_{t}\mathrm{d}t+K^{-}V_{t}\mathrm{d}\Lambda _{-}=\left(
J-I\right) V_{t}\mathrm{d}\Lambda +L_{+}V_{t}\mathrm{d}\Lambda ^{+}
\label{1.6}
\end{equation}%
where $\Lambda ^{+}\left( t\right) $ is the creation process, corresponding
to the annihilation $\Lambda _{-}\left( t\right) $ on the interval $[0,t)$,
and $\Lambda \left( t\right) $ is the number of quanta on this interval.
These canonical quantum stochastic processes, representing the quantum noise
with respect to the vacuum state $|0\rangle $ of the Fock space $\mathcal{F}$
over the single-quantum Hilbert space $L^{2}\left( \mathbb{R}_{+}\right) $
of square-integrable functions of $t\in \lbrack 0,\infty )$, are formally
given in \cite{14} by the integrals 
\begin{equation*}
\Lambda _{-}\left( t\right) =\int_{0}^{t}\Lambda _{-}^{r}\mathrm{d}r,\quad
\Lambda ^{+}\left( t\right) =\int_{0}^{t}\Lambda _{r}^{+}\mathrm{d}r,\quad
\Lambda \left( t\right) =\int_{0}^{t}\Lambda _{r}^{+}\Lambda _{-}^{r}\mathrm{%
d}r,
\end{equation*}%
where $\Lambda _{-}^{r},\Lambda _{r}^{+}$ are the generalized quantum
one-dimensional fields in $\mathcal{F}$, satisfying the canonical
commutation relations 
\begin{equation*}
\left[ \Lambda _{-}^{r},\Lambda _{s}^{+}\right] =\delta \left( s-r\right)
I,\quad \left[ \Lambda _{-}^{r},\Lambda _{-}^{s}\right] =0=\left[ \Lambda
_{r}^{+},\Lambda _{s}^{+}\right] .
\end{equation*}%
They can be defined by the independent increments with 
\begin{equation}
\langle 0|\mathrm{d}\Lambda _{-}|0\rangle =0,\quad \langle 0|\mathrm{d}%
\Lambda ^{+}|0\rangle =0,\quad \langle 0|\mathrm{d}\Lambda |0\rangle =0
\label{1.7}
\end{equation}%
and the noncommutative multiplication table 
\begin{equation}
\mathrm{d}\Lambda \mathrm{d}\Lambda =\mathrm{d}\Lambda ,\quad \mathrm{d}%
\Lambda _{-}\mathrm{d}\Lambda =\mathrm{d}\Lambda _{-},\quad \mathrm{d}%
\Lambda \mathrm{d}\Lambda ^{+}=\mathrm{d}\Lambda ^{+},\quad \mathrm{d}%
\Lambda _{-}\mathrm{d}\Lambda ^{+}=\mathrm{d}tI  \label{1.8}
\end{equation}%
with all other products being zero: $\mathrm{d}\Lambda \mathrm{d}\Lambda
_{-}=\mathrm{d}\Lambda ^{+}\mathrm{d}\Lambda =\mathrm{d}\Lambda ^{+}\mathrm{d%
}\Lambda _{-}=0$. The standard Poisson process $\mathrm{P}$ as well as the
Wiener process $\mathrm{Q}$ can be represented in $\mathfrak{F}$ by the
linear combinations \cite{16} 
\begin{equation}
\mathrm{P}\left( t\right) =\Lambda \left( t\right) +i\left( \Lambda
^{+}\left( t\right) -\Lambda _{-}\left( t\right) \right) ,\quad \mathrm{Q}%
\left( t\right) =\Lambda ^{+}\left( t\right) +\Lambda _{-}\left( t\right) ,
\label{1.9}
\end{equation}%
so the equation (\ref{1.6}) corresponds to the stochastic diffusion equation
(\ref{1.3}) if $J=I$, $L_{+}=L=-K^{-}$, and it corresponds to the stochastic
jump equation (\ref{1.1}) if $J=I+L$, $L_{+}=iL=K^{-}$. The quantum
stochastic equation for $\phi _{t}\left( B\right) =V_{t}^{\dagger }BV_{t}$
has the following general form 
\begin{equation*}
\mathrm{d}\phi _{t}\left( B\right) +\phi _{t}\left( K^{\dagger
}B+BK-L^{-}BL_{+}\right) \mathrm{d}t=\phi _{t}\left( J^{\dagger }BJ-B\right) 
\mathrm{d}\Lambda
\end{equation*}%
\begin{equation}
+\phi _{t}\left( J^{\dagger }BL_{+}-K_{+}B\right) +\phi _{t}\left(
L^{-}BJ-BK^{-}\right) \mathrm{d}\Lambda _{-},  \label{1.10}
\end{equation}%
where $L^{-}=L_{+}^{\dagger },K_{+}^{\dagger }=K^{-}$, coinciding with
either (\ref{1.5}) or with (\ref{1.4}) in the particular cases. The equation
(\ref{1.10}) is obtained from (\ref{1.6}) by using the It\^{o} formula (\ref%
{1.2}) with the multiplication table (\ref{1.8}). The sub-filtering
condition $K+K^{\dagger }\leq L^{-}L_{+}$ for the equation (\ref{1.6})
defines in both cases the positive operator-valued process $R_{t}=\phi
_{t}\left( I\right) $ as a sub-martingale with $R_{0}=I$, or a martingale in
the case $K+K^{\dagger }=L^{-}L_{+}$. In the particular case 
\begin{equation*}
J=S,\quad K^{-}=L^{-}S,\quad L_{+}=SK_{+},\quad S^{\dagger }S=I,
\end{equation*}%
corresponding to the Hudson--Evans flow if $S^{\dagger }=S^{-1}$, the
evolution is isometric, and identity preserving, $\phi _{t}\left( I\right)
=I $ in the case of bounded $K$ and $L$.

In the next sections we define a multidimensional analog of the quantum
stochastic equation (\ref{1.10}) and will show that the suggested general
structure of its generator indeed follows just from the property of complete
positivity of the map $\phi _t$ for all $t>0$ and the normalization
condition $\phi _t\left( I\right) =M_t$ to a form-valued sub-martingale with
respect to the natural filtration of the quantum noise in the Fock space $%
\mathfrak{F}$ .

\section{The Generators of Quantum Filtering Cocycles.}

The quantum filtering dynamics over an operator algebra $\mathcal{B}%
\subseteq \mathcal{B}\left( \mathcal{H}\right) $ is described by a one
parameter cocycles: $\phi =\left( \phi _t\right) _{t>0}$ of linear
completely positive stochastic maps $\phi _t\left( \omega \right) :\mathcal{B%
}\rightarrow \mathcal{B}$. The cocycle condition 
\begin{equation}
\phi _s\left( \omega \right) \circ \phi _r\left( \omega ^s\right) =\phi
_{r+s}\left( \omega \right) ,\quad \forall r,s>0  \label{2.1}
\end{equation}
means the stationarity, with respect to the shift $\omega ^s=\left\{ \omega
\left( t+s\right) \right\} $ of a given stochastic process $\omega =\left\{
\omega \left( t\right) \right\} $. Such maps are in general unbounded, but
normalized, $\phi _t\left( I\right) =M_t$ to an operator-valued martingale $%
M_t=\epsilon _t\left[ M_s\right] \geq 0$ with $M_0=1$, or a positive
submartingale: $M_t\geq \epsilon _t\left[ M_s\right] $, for all $s>t$, .

Now we give a noncommutative generalization of the quantum stochastic CP
cocycles, which was suggested in \cite{15} even for the nonlinear case. The
stochastically differentiable family $\phi $ with respect to a quantum
stationary process, with independent increments $\Lambda ^s\left( t\right)
=\Lambda \left( t+s\right) -\Lambda \left( s\right) $ generated by a finite
dimensional It\^{o} algebra is described by the quantum stochastic equation 
\begin{equation}
\mathrm{d}\phi _t\left( Y\right) =\phi _t\circ \lambda _\nu ^\mu \left(
Y\right) \mathrm{d}\Lambda _\mu ^\nu :=\sum_{\mu ,\nu }\phi _t\left( \lambda
_\nu ^\mu \left( Y\right) \right) \mathrm{d}\Lambda _\mu ^\nu ,\qquad \text{ 
}Y\in \mathcal{B}\qquad  \label{2.2}
\end{equation}
with the initial condition $\phi _0\left( Y\right) =Y$, for all $Y\in 
\mathcal{B}$. Here $\Lambda _\mu ^\nu \left( t\right) $ with $\mu \in
\left\{ -,1,...,d\right\} $, $\nu \in \left\{ +,1,...,d\right\} $ are the
standard time $\Lambda _{-}^{+}\left( t\right) =tI$, annihilation $\Lambda
_{-}^m\left( t\right) $, creation $\Lambda _n^{+}\left( t\right) $ and
exchange-number $\Lambda _n^m\left( t\right) =\mathrm{N}_n^m\left( t\right) $
operator integrators with $m,n\in \left\{ 1,...,d\right\} $. The
infinitesimal increments $\mathrm{d}\Lambda _\nu ^\mu \left( t\right)
=\Lambda _\nu ^{t\mu }\left( \mathrm{d}t\right) $ are formally defined by
the Hudson-Parthasarathy multiplication table \cite{16} and the $\flat $
-property \cite{3}, 
\begin{equation}
\mathrm{d}\Lambda _\mu ^\beta \mathrm{d}\Lambda _\gamma ^\nu =\delta _\gamma
^\beta \mathrm{d}\Lambda _\mu ^\nu ,\qquad \text{ }\Lambda ^{\flat }=\Lambda
,\qquad  \label{2.3}
\end{equation}
where $\delta _\gamma ^\beta $ is the usual Kronecker delta restricted to
the indices $\beta \in \left\{ -,1,...,d\right\} ,\quad \gamma \in \left\{
+,1,...,d\right\} $ and $\Lambda _{-\nu }^{\flat \mu }=\Lambda _{-\mu }^{\nu
\dagger }$ with respect to the reflection $-(-)=+,$ $-(+)=-$ of the indices $%
\left( -,+\right) $ only. The linear maps $\lambda _\nu ^\mu :\mathcal{B}%
\rightarrow \mathcal{B}$ for the $*$ -cocycles $\phi _t^{*}=\phi _t$, where $%
\phi _t^{*}\left( Y\right) =\phi _t\left( Y^{\dagger }\right) ^{\dagger }$,
should obviously satisfy the $\flat $ -property $\lambda ^{\flat }=\lambda $%
, where $\lambda _{-\mu }^{\flat \nu }=\lambda _{-\nu }^{\mu *}$, $\lambda
_\nu ^{\mu *}\left( Y\right) =\lambda _\nu ^\mu \left( Y^{\dagger }\right)
^{\dagger }$. If the coefficients $b_\nu ^\mu =\lambda _\nu ^\mu \left(
Y\right) $ are independent of $t$, $\phi $ satisfies the cocycle property $%
\phi _s\circ \phi _r^s=\phi _{s+r}$, where $\phi _t^s$ is the solution to (%
\ref{2.2}) with $\Lambda _\nu ^\mu \left( t\right) $ replaced by $\Lambda
_\nu ^{s\mu }\left( t\right) $. Define the $\left( d+2\right) \times \left(
d+2\right) $ matrix $\mathbf{a}=\left[ a_\nu ^\mu \right] $ also for $\mu =+$
and $\nu =-$, by 
\begin{equation*}
\lambda _\nu ^{+}\left( Y\right) =0=\lambda _{-}^\mu \left( Y\right) ,\qquad
\forall Y\in \mathcal{B},
\end{equation*}
and then one can extend the summation in (\ref{2.2}) so it is also over $\mu
=+$, and $\nu =-$. By such an extension the multiplication table for $%
\mathrm{d}\Lambda \left( \mathbf{a}\right) =a_\nu ^\mu \mathrm{d}\Lambda
_\mu ^\nu $ can be written as 
\begin{equation}
\mathrm{d}\Lambda \left( \mathbf{a}\right) ^{\dagger }\mathrm{d}\Lambda
\left( \mathbf{a}\right) =\mathrm{d}\Lambda \left( \mathbf{a}^{\flat }%
\mathbf{a}\right)  \label{2.4}
\end{equation}
in terms of the usual matrix product $\left( \mathbf{ba}\right) _\nu ^\mu
=b_\lambda ^\mu a_\nu ^\lambda $ and the involution $\mathbf{a\mapsto a}%
^{\flat }=\mathbf{b},\mathbf{b}^{\flat }=\mathbf{a}$ can be obtained by the
pseudo-Hermitian conjugation $a_\beta ^{\flat \nu }=g_{\beta \mu }a_\gamma
^{\mu \dagger }g^{\gamma \nu }$ respectively to the indefinite Minkowski
metric tensor $\mathrm{g}=\left[ g_{\mu \nu }\right] $ and its inverse $%
\mathrm{g}^{-1}=\left[ g^{\mu \nu }\right] $, given by $g^{\mu \nu }=\delta
_{-\nu }^\mu I=g_{\mu \nu }$.

Let us prove that the "spatial" part $\boldsymbol{\gamma }=\left( \gamma
_{\nu }^{\mu }\right) _{\nu \neq -}^{\mu \neq +}$ of $\gamma =\lambda
+\delta $, called the quantum stochastic germ for the representation $\delta
:B\mapsto \left( B\delta _{\nu }^{\mu }\right) _{\nu \neq -}^{\mu \neq +}$,
must be completely stochastically dissipative for a CP cocycle $\phi $ in
the following sense.

\begin{theorem}
Suppose that the quantum stochastic equation (\ref{2.2}) with $\phi _0\left(
B\right) =B$ has a CP solution $\phi _t,t>0$. Then the germ-map $\boldsymbol{%
\gamma }=\left( \lambda _\nu ^\mu +\delta _\nu ^\mu \right) _{\nu =+,\bullet
}^{\mu =-,\bullet }$ is conditionally completely positive 
\begin{equation*}
\sum_k\boldsymbol{\iota }\left( B_k\right) \boldsymbol{\eta }_k=0\Rightarrow
\sum_{k,l}\langle \boldsymbol{\eta }_k|\boldsymbol{\gamma }\left(
B_k^{\dagger }B_l\right) \boldsymbol{\eta }_l\rangle \geq 0
\end{equation*}
Here $\boldsymbol{\eta }\in \mathcal{H}\oplus \mathcal{H}^{\bullet },%
\mathcal{H}^{\bullet }=\mathcal{H}\otimes \mathbb{C}^d$, and $\boldsymbol{%
\iota }=\left( \iota _\nu ^\mu \right) _{\nu =+,\bullet }^{\mu =-,\bullet }$
is the degenerate representation $\iota _\nu ^\mu \left( B\right) =B\delta
_\nu ^{+}\delta _{-}^\mu $, written both with $\boldsymbol{\gamma }$ in the
matrix form as 
\begin{equation}
\boldsymbol{\gamma }=\left( 
\begin{array}{cc}
\gamma & \gamma _{\bullet } \\ 
\gamma ^{\bullet } & \gamma _{\bullet }^{\bullet }%
\end{array}
\right) ,\qquad \text{ }\boldsymbol{\iota }\left( B\right) =\left( 
\begin{array}{cc}
B & 0 \\ 
0 & 0%
\end{array}
\right) ,\qquad  \label{2.5}
\end{equation}
where $\gamma =\lambda _{+}^{-},\quad $ $\gamma ^m=\lambda _{+}^m,\quad $ $%
\gamma _n=\lambda _n^{-},\quad \gamma _n^m=\delta _n^m+\lambda _n^m$ with $%
\delta _n^m\left( B\right) =B\delta _n^m$ such that 
\begin{equation}
\gamma \left( B^{\dagger }\right) =\gamma \left( B\right) ^{\dagger },\qquad 
\text{ }\gamma ^n\left( B^{\dagger }\right) =\gamma _n\left( B\right)
^{\dagger },\qquad \text{ }\gamma _n^m\left( B^{\dagger }\right) =\gamma
_m^n\left( B\right) ^{\dagger }  \label{2.6}
\end{equation}
%TCIMACRO{\TeXButton{Proof}{\proof} }%
%BeginExpansion
\proof
%EndExpansion
Let us denote by $\mathcal{D}$ the $\mathcal{H}$-span $\left\{ \sum_f\xi
^f\otimes f^{\otimes }\left| \xi ^f\in \mathcal{H},f^{\bullet }\in \mathbb{C}%
^d\otimes L^2\left( \mathbb{R}_{+}\right) \right. \right\} $ of coherent
(exponential) functions $f^{\otimes }\left( \tau \right) =\bigotimes_{t\in
\tau }f^{\bullet }\left( t\right) $, given for each finite subset $\tau
=\left\{ t_1,...,t_n\right\} \subseteq \mathbb{R}_{+}$ by tensor products $%
f^{n_1,...,n_N}\left( \tau \right) =f^{n_1}\left( t_1\right)
...f^{n_N}\left( t_N\right) $, where $f^n,n=1,...,d$ are square-integrable
complex functions on $\mathbb{R}_{+}$ and $\xi ^f=0$ for almost all $%
f^{\bullet }=\left( f^n\right) $. The co-isometric shift $T_s$ intertwining $%
\text{A}^s\left( t\right) $ with $\text{A}\left( t\right) =T_s\text{A}%
^s\left( t\right) T_s^{\dagger }$ is defined on $\mathcal{D}$ by $T_s\left(
\eta \otimes f^{\otimes }\right) \left( \tau \right) =\eta \otimes
f^{\otimes }\left( \tau +s\right) $. The complete positivity of the quantum
stochastic adapted map $\phi _t$ into the $\mathcal{D}$-forms $\left\langle
\chi \right| \left. \phi _t\left( B\right) \psi \right\rangle $, for $\chi
,\psi \in \mathcal{D}$ can be obviously written as 
\begin{equation}
\sum_{X,Z}\sum_{f,h}\left\langle \xi _X^f\right| \left. \phi _t\left(
f^{\bullet },X^{\dagger }Z,h^{\bullet }\right) \xi _Z^h\right\rangle \geq
0,\qquad  \label{2.7}
\end{equation}
where 
\begin{equation*}
\left\langle \eta \right| \left. \phi _t\left( f^{\bullet },B,h^{\bullet
}\right) \eta \right\rangle =\left\langle \eta \otimes f^{\otimes }\right|
\left. \phi _t\left( B\right) \eta \otimes h^{\otimes }\right\rangle
e^{-\int_t^\infty f^{\bullet }\left( s\right) ^{\dagger }h^{\bullet }\left(
s\right) \mathrm{d}s},
\end{equation*}
$\xi _B^f\neq 0$ for a finite sequence of $B_k\in \mathcal{B}$, and for a
finite sequence of $f_l^{\bullet }=\left( f_l^1,...,f_l^d\right) $. If the $%
\mathcal{D}$-form $\phi _t\left( B\right) $ satisfies the stochastic
equation (\ref{2.2}), the $\mathcal{H}$-form $\phi _t\left( f^{\bullet
},B,h^{\bullet }\right) $ satisfies the differential equation 
\begin{equation}
\frac{\mathrm{d}}{\mathrm{d}t}\phi _t\left( f^{\bullet },B,h^{\bullet
}\right) =f^{\bullet }\left( t\right) ^{\dagger }h^{\bullet }\left( t\right)
\phi _t\left( f^{\bullet },B,h^{\bullet }\right) +\phi _t\left( f^{\bullet
},\lambda _{+}^{-}\left( B\right) ,h^{\bullet }\right)  \label{2.8}
\end{equation}
\begin{eqnarray*}
&&+\sum_{m=1}^df^m\left( t\right) ^{*}\phi _t\left( f^{\bullet },\lambda
_{+}^m\left( B\right) ,h^{\bullet }\right) +\sum_{n=1}^dh^n\left( t\right)
\phi _t\left( f^{\bullet },\lambda _n^{-}\left( B\right) ,h^{\bullet }\right)
\\
&&+\sum_{m,n=1}^df^m\left( t\right) ^{*}h^n\left( t\right) \phi _t\left(
f^{\bullet },\lambda _n^m\left( B\right) ,h^{\bullet }\right) ,
\end{eqnarray*}
where $f^{\bullet }\left( t\right) ^{\dagger }h^{\bullet }\left( t\right)
=\sum_{n=1}^df^n\left( t\right) ^{*}h^n\left( t\right) $. The positive
definiteness, (\ref{2.7}), ensures the conditional positivity 
\begin{equation}
\sum_f\sum_BB\xi _B^f=0\Rightarrow \sum_{X,Z}\sum_{f,h}\left\langle \xi
_X^f\right| \left. \gamma \left( f^{\bullet },X^{\dagger }Z,h^{\bullet
}\right) \xi _Z^h\right\rangle \geq 0  \label{2.9}
\end{equation}
of the form $\gamma _t\left( f^{\bullet },B,h^{\bullet }\right) =\frac
1t\left( \phi _t\left( f^{\bullet },B,h^{\bullet }\right) -B\right) $ for
each $t>0$ and of the limit $\gamma _0$ at $t\downarrow 0$, coinciding with
the quadratic form 
\begin{equation}
\left. \frac{\mathrm{d}}{\mathrm{d}t}\phi _t\left( f^{\bullet },B,h^{\bullet
}\right) \right| _{t=0}=\sum_{m,n}\bar{a}^m\gamma _n^m\left( B\right)
c^n+\sum_m\bar{a}^m\gamma ^m\left( B\right) +\sum_n\gamma _n\left( B\right)
c^n+\gamma \left( B\right) ,  \label{2.10}
\end{equation}
where $a^{\bullet }=f^{\bullet }\left( 0\right) ,\quad c^{\bullet
}=h^{\bullet }\left( 0\right) $, and the $\gamma $'s are defined in (\ref%
{2.5}). Hence the form 
\begin{equation*}
\sum_{X,Z}\sum_{\mu ,\nu }\left\langle \eta _X^\mu \right| \gamma _\nu ^\mu
\left( X^{\dagger }Z\right) \eta _Z^\nu \rangle
:=\sum_{X,Z}\sum_{m,n}\left\langle \eta _X^m\right| \left. \gamma _n^m\left(
X^{\dagger }Z\right) \eta _Z^n\right\rangle
\end{equation*}
\begin{equation*}
+\sum_{X,Z}\left( \sum_n\left\langle \eta _X\right| \left. \gamma _n\left(
X^{\dagger }Z\right) \eta _Z^n\right\rangle +\sum_m\left\langle \eta
_X^m\right| \left. \gamma ^m\left( X^{\dagger }Z\right) \eta _Z\right\rangle
+\left\langle \eta _X\right| \gamma \left( X^{\dagger }Z\right) \left| \eta
_Z\right\rangle \right)
\end{equation*}
with $\eta =\sum_f\xi ^f,\quad \eta ^{\bullet }=\sum_f\xi ^f\otimes
a_f^{\bullet }$, where $a_f^{\bullet }=f^{\bullet }\left( 0\right) $, is
positive if $\sum_BB\eta _B=0$. The components $\eta $ and $\eta ^{\bullet }$
of these vectors are independent because for any $\eta \in \mathcal{H}$ and $%
\eta ^{\bullet }=\left( \eta ^1,...,\eta ^d\right) \in \mathcal{H}\otimes 
\mathbb{C}^d$ there exists such a function $a^{\bullet }\mapsto \xi ^a$ on $%
\mathbb{C}^d$ with a finite support, that $\sum_a\xi ^a=\eta ,\quad
\sum_a\xi ^a\otimes a^{\bullet }=\eta ^{\bullet }$, namely, $\xi ^a=0$ for
all $a^{\bullet }\in \mathbb{C}^d$ except $a^{\bullet }=0$, for which $\xi
^a=\eta -\sum_{n=1}^d\eta ^n$ and $a^{\bullet }=e_n^{\bullet }$, the $n$-th
basis element in $\mathbb{C}^d$, for which $\xi ^a=\eta ^n$. This proves the
complete positivity of the matrix form $\boldsymbol{\gamma }$, with respect
to the matrix representation $\boldsymbol{\iota }$ defined in (\ref{2.5}) on
the ket-vectors $\boldsymbol{\eta }=\left( \eta ^\mu \right) $. 
%TCIMACRO{\TeXButton{End Proof}{\endproof}}%
%BeginExpansion
\endproof%
%EndExpansion
\end{theorem}

\section{A Dilation Theorem for the Form-Generator.}

The conditional positivity of the structural map $\boldsymbol{\gamma }$ with
respect to the degenerate representation $\boldsymbol{\iota }$ written in
the matrix form (\ref{2.6}) obviously implies the positivity of the
dissipation form 
\begin{equation}
\sum_{X,Z}\left\langle \boldsymbol{\eta }_X\right| \boldsymbol{\Delta }%
\left( X,Z\right) \left. \boldsymbol{\eta }_Z\right\rangle
:=\sum_{k,l}\sum_{\mu ,\nu }\left\langle \eta _k^\mu \right| \Delta _\nu
^\mu \left( B_k,B_l\right) \left. \eta _l^\nu \right\rangle ,\qquad
\label{3.1}
\end{equation}
where $\eta ^{-}=\eta =\eta ^{+}$ and $\eta _k=\eta _{B_k}$ for any (finite)
sequence $B_k\in \mathcal{B}$, $k=1,2,...$, corresponding to non-zero $%
\boldsymbol{\eta }_B=\eta _B\oplus \eta _B^{\bullet },\eta _B\in \mathcal{H}%
,\eta _B^{\bullet }\in \mathcal{H}^{\bullet }$. Here $\boldsymbol{\Delta }%
=\left( \Delta _\nu ^\mu \right) _{\nu =+,\bullet }^{\mu =-,\bullet }$ is
the dissipator matrix, 
\begin{equation*}
\boldsymbol{\Delta }\left( X,Z\right) =\boldsymbol{\gamma }\left( X^{\dagger
}Z\right) -\boldsymbol{\iota }\left( X\right) ^{\dagger }\boldsymbol{\gamma }%
\left( Z\right) -\boldsymbol{\gamma }\left( X\right) ^{\dagger }\boldsymbol{%
\iota }\left( Z\right) +\boldsymbol{\iota }\left( X\right) ^{\dagger }%
\boldsymbol{\gamma }\left( I\right) \boldsymbol{\iota }\left( Z\right) ,
\end{equation*}
given by the elements 
\begin{eqnarray}
\Delta _n^m\left( X,Z\right) &=&\lambda _n^m\left( X^{\dagger }Z\right)
+X^{\dagger }Z\delta _n^m,  \label{3.2} \\
\Delta _n^{-}\left( X,Z\right) &=&\lambda _n^{-}\left( X^{\dagger }Z\right)
-X^{\dagger }\lambda _n^{-}\left( Z\right) =\Delta _{+}^n\left( Z,X\right)
^{\dagger }  \notag \\
\Delta _{+}^{-}\left( X,Z\right) &=&\lambda _{+}^{-}\left( X^{\dagger
}Z\right) -X^{\dagger }\lambda _{+}^{-}\left( Z\right) -\lambda
_{+}^{-}\left( X^{\dagger }\right) Z+X^{\dagger }DZ,  \notag
\end{eqnarray}
where $D=\lambda _{+}^{-}\left( I\right) \leq 0$ ($D=0$ for the case of the
martingale $M_t$ ). This means that the matrix-valued map $\gamma _{\bullet
}^{\bullet }=\left[ \gamma _n^m\right] $, is completely positive, and as
follows from the next theorem, at least for the algebra $\mathcal{B}=%
\mathcal{B}\left( \mathcal{H}\right) $ the maps $\gamma $, $\gamma ^m$, $%
\gamma _n$ have the following form 
\begin{eqnarray}
\gamma ^m\left( B\right) &=&\varphi ^m\left( B\right) -K_m^{\dagger
}B,\qquad \text{ }\gamma _n\left( B\right) =\varphi _n\left( B\right)
-BK_n\qquad  \label{3.3} \\
\gamma \left( B\right) &=&\varphi \left( B\right) -K^{\dagger }B-BK,\qquad 
\text{ }\varphi \left( I\right) \leq K+K^{\dagger }  \notag
\end{eqnarray}
where $\boldsymbol{\varphi }=\left( \varphi _\nu ^\mu \right) _{\nu \neq
-}^{\mu \neq +}$ is a completely positive bounded map from $\mathcal{B}$
into the matrices of operators with the elements $\varphi _n^m=\gamma
_n^m,\varphi _{+}^m=\varphi ^m,\varphi _n^{-}=\varphi _n,\varphi
_{+}^{-}=\varphi :\mathcal{B}\rightarrow \mathcal{B}$.

In order to make the formulation of the dilation theorem as concise as
possible, we need the notion of the $\flat $-representation of the algebra $%
\mathcal{B}$ in the operator algebra $\mathcal{A}\left( \mathcal{E}\right) $
of a pseudo-Hilbert space $\mathcal{E}=\mathcal{H}\oplus \mathcal{H}^{\circ
}\oplus \mathcal{H}$ with respect to the indefinite metric 
\begin{equation}
\left( \xi \right| \left. \xi \right) =2\mathrm{Re}\left( \xi ^{-}\right|
\left. \xi ^{+}\right) +\left\| \xi ^{\circ }\right\| ^2+\left\| \xi
^{+}\right\| _D^2  \label{3.4}
\end{equation}
for the triples $\xi =\left( \xi ^\mu \right) ^{\mu =-,\circ ,+}\in \mathcal{%
E}$, where $\xi ^{-},\xi ^{+}\in \mathcal{H},\quad \xi ^{\circ }\in \mathcal{%
H}^{\circ },\quad \mathcal{H}^{\circ }$ is a pre-Hilbert space, and $\left\|
\eta \right\| _D^2=\left\langle \eta \right| \left. D\eta \right\rangle $.
The operators $A\in \mathcal{A}\left( \mathcal{E}\right) $ are given by $%
3\times 3$-block-matrices $\left[ A_\nu ^\mu \right] _{\nu =-,\circ ,+}^{\mu
=-,\circ ,+}$, having the Pseudo-Hermitian adjoints $\left( \xi |A^{\flat
}\xi \right) =\left( A\xi |\xi \right) $, which are defined by the Hermitian
adjoints $A_\nu ^{\dagger \mu }=A_\mu ^{\nu \dagger }$as $A^{\flat
}=G^{-1}A^{\dagger }G$ respectively to the indefinite metric tensor $G=\left[
G_{\mu \nu }\right] $ and its inverse $G^{-1}=\left[ G^{\mu \nu }\right] $,
given by 
\begin{equation}
\mathbf{G}=\left[ 
\begin{array}{ccc}
0 & 0 & I \\ 
0 & I_{\circ }^{\circ } & 0 \\ 
I & 0 & D%
\end{array}
\right] ,\qquad \mathbf{G}^{-1}=\left[ 
\begin{array}{ccc}
-D & 0 & I \\ 
0 & I_{\circ }^{\circ } & 0 \\ 
I & 0 & 0%
\end{array}
\right]  \label{3.5}
\end{equation}
with an arbitrary $D$, where $I_{\circ }^{\circ }$ is the identity operator
in $\mathcal{H}^{\circ }$, being equal $I_{\bullet }^{\bullet }=\left[
I\delta _n^m\right] _{n=1,...,d}^{m=1,...,d}$ in the case of $\mathcal{H}%
^{\circ }=\mathcal{H}\otimes \mathbb{C}^d=\mathcal{H}^{\bullet }$.

\begin{theorem}
The following are equivalent:

\begin{enumerate}
\item[(i)] The dissipation form (\ref{3.1}), defined by the $\flat $-map $%
\lambda $ with $\lambda _{+}^{-}\left( I\right) =D$, is positive definite: $%
\sum_{X,Z}\left\langle \boldsymbol{\eta }_X\right| \boldsymbol{\Delta }%
\left( X,Z\right) \left. \boldsymbol{\eta }_Z\right\rangle \geq 0$.

\item[(ii)] There exists a pre-Hilbert space $\mathcal{H}^{\circ }$, a
unital $\dagger $- representation $j$ of $\mathcal{B}$ in $\mathcal{B}\left( 
\mathcal{H}^{\circ }\right) $, 
\begin{equation}
j\left( B^{\dagger }B\right) =j\left( B\right) ^{\dagger }j\left( B\right)
,\quad j\left( I\right) =I,  \label{3.6}
\end{equation}
a $\left( j,i\right) $-derivation of $\mathcal{B}$ with $i\left( B\right) =B$%
, 
\begin{equation}
k\left( B^{\dagger }B\right) =j\left( B\right) ^{\dagger }k\left( B\right)
+k\left( B^{\dagger }\right) B,  \label{3.7}
\end{equation}
having values in the operators $\mathcal{H}\rightarrow \mathcal{H}^{\circ }$%
, the adjoint map $k^{*}\left( B\right) =k\left( B^{\dagger }\right)
^{\dagger }$, with the property 
\begin{equation*}
k^{*}\left( B^{\dagger }B\right) =B^{\dagger }k^{*}\left( B\right)
+k^{*}\left( B^{\dagger }\right) j\left( B\right)
\end{equation*}
of $\left( i,j\right) $-derivation in the operators $\mathcal{H}^{\circ
}\rightarrow \mathcal{H}$, and a map $l:\mathcal{B}\rightarrow \mathcal{B}$
having the coboundary property 
\begin{equation}
l\left( B^{\dagger }B\right) =B^{\dagger }l\left( B\right) +l\left(
B^{\dagger }\right) B+k^{*}\left( B^{\dagger }\right) k\left( B\right) ,
\label{3.8}
\end{equation}
with the adjoint $l^{*}\left( B\right) =l\left( B\right) +\left[ D,B\right] $%
, such that $\gamma \left( B\right) =l\left( B\right) +DB$, 
\begin{equation*}
\gamma _n\left( B^{\dagger }\right) =k\left( B\right) ^{\dagger }L_n^{\circ
}+B^{\dagger }L_n^{-}=\gamma ^n\left( B\right) ^{\dagger },
\end{equation*}
and $\gamma _n^m\left( B\right) =L_m^{\circ \dagger }j\left( B\right)
L_n^{\circ }$ for some operators $L_n^{\circ }:\mathcal{H}\rightarrow 
\mathcal{H}^{\circ }$ having the adjoints $L_n^{\circ \dagger }$ on $%
\mathcal{H}^{\circ }$ and $L_n^{-}\in \mathcal{B}$.

\item[(iii)] There exists a pseudo-Hilbert space, $\mathcal{E}$, a unital $%
\flat $-representation $\jmath :\mathcal{B\rightarrow A}\left( \mathcal{E}%
\right) $, and a linear operator $\mathbf{L}:\mathcal{H}\oplus \mathcal{H}%
^{\bullet }\rightarrow \mathcal{E}$ such that 
\begin{equation}
\mathbf{L}^{\flat }\jmath \left( B\right) \mathbf{L}=\boldsymbol{\gamma }%
\left( B\right) ,\qquad \forall B\in \mathcal{B}.  \label{3.10}
\end{equation}

\item[(iv)] The structural map $\boldsymbol{\gamma }=\boldsymbol{\lambda }+%
\boldsymbol{\delta }$ is conditionally completely positive with respect to
the matrix representation $\boldsymbol{\iota }$ in (\ref{2.5}).
\end{enumerate}
\end{theorem}

%TCIMACRO{\TeXButton{Proof}{\proof} }%
%BeginExpansion
\proof
%EndExpansion
The implication (i)$\Rightarrow $(ii) generalizes the Evans-Lewis Theorem 
\cite{17}, and its proof is similar to the proof of the dilation theorem in 
\cite{18}. Let $\mathcal{H}^{\circ }$ be the pre-Hilbert space of Kolmogorov
decomposition $\boldsymbol{\Delta }\left( X,Z\right) =\boldsymbol{k}\left(
X\right) ^{\dagger }\boldsymbol{k}\left( Z\right) $. It is defined as the
quotient space $\mathcal{H}^{\circ }=\mathcal{K}/\mathcal{I}$ of the $%
\mathcal{H}$-span $\mathcal{K}=\left\{ \left( \boldsymbol{\eta }_B\right)
_{B\in B}\right\} $, where $\boldsymbol{\eta }_B\in \mathcal{H}\oplus 
\mathcal{H}^{\bullet }$ is not equal zero only for a finite number of $B\in
B $, with respect to the kernel 
\begin{equation*}
\mathcal{I}=\left\{ \left( \boldsymbol{\eta }_B\right) _{B\in B}\in \mathcal{%
K}|\sum_{X,Z}\left\langle \boldsymbol{\eta }_X\right| \boldsymbol{\Delta }%
\left( X,Z\right) \left. \boldsymbol{\eta }_Z\right\rangle =0\right\}
\end{equation*}
of the positive-definite form (\ref{3.1}). The operators $\boldsymbol{k}%
\left( B\right) ^{\dagger }:\mathcal{H}^{\circ }\rightarrow \mathcal{H}%
\oplus \mathcal{H}^{\bullet }$ are defined on the classes $\eta ^{\circ }$
of $\left( \boldsymbol{\eta }_X\right) _{X\in B}\in \mathcal{K}$ as the
adjoint 
\begin{equation*}
\left\langle \boldsymbol{k}\left( B\right) ^{\dagger }\eta ^{\circ }|%
\boldsymbol{\eta }\right\rangle =\sum_X\left\langle \boldsymbol{\eta }_X|%
\boldsymbol{\Delta }\left( X,B\right) \boldsymbol{\eta }\right\rangle
\end{equation*}
to the bounded operators $\boldsymbol{k}\left( B\right) :\mathcal{H}\oplus 
\mathcal{H}^{\bullet }\rightarrow \mathcal{H}^{\circ }$, mapping the pairs $%
\boldsymbol{\eta }=\eta \oplus \eta ^{\bullet }$ into the equivalence
classes $\eta ^{\circ }\left( B\right) =k\left( B\right) \eta +k_{\bullet
}\left( B\right) \eta ^{\bullet }$ of $\left( \delta _Z\left( B\right) 
\boldsymbol{\eta }\right) _{Z\in B}$, where $\delta _Z\left( B\right) =1$ if 
$B=Z$, otherwise $\delta _Z\left( B\right) =0$. Let us define a linear
operator $j\left( B\right) $ on $\mathcal{H}^{\circ }$ by 
\begin{equation*}
j\left( B\right) \sum_Z\left( k\left( Z\right) \eta +k_{\bullet }\left(
Z\right) \eta ^{\bullet }\right) =\sum_Z\left( k\left( BZ\right) \eta
-k\left( B\right) Z\eta +k_{\bullet }\left( BZ\right) \eta ^{\bullet
}\right) .
\end{equation*}
Obviously $j\left( XB\right) =j\left( X\right) j\left( B\right) $, $j\left(
I\right) =I$ because $k\left( I\right) =0$ and as follows from the
definition of the dissipation form, $j\left( B\right) ^{\dagger }=j\left(
B^{\dagger }\right) $ for all $B\in \mathcal{B}$. Thus $j$ is a unital $%
\dagger $-representation, $k$ is a $\left( j,i\right) $-cocycle, and $%
k_{\bullet }\left( B\right) =j\left( B\right) L_{\bullet }^{\circ }$, where $%
L_{\bullet }^{\circ }=k_{\bullet }\left( I\right) $. Moreover, as 
\begin{eqnarray*}
\gamma \left( B^{\dagger }B\right) +B^{\dagger }\gamma \left( I\right) B
&=&B^{\dagger }\gamma \left( B\right) +\gamma \left( B^{\dagger }\right)
B+k\left( B\right) ^{\dagger }k\left( B\right) , \\
\gamma _{\bullet }^{\bullet }\left( B^{\dagger }B\right) &=&k_{\bullet
}\left( B\right) ^{\dagger }k_{\bullet }\left( B\right) , \\
\gamma _{\bullet }\left( B^{\dagger }B\right) -B^{\dagger }\gamma _{\bullet
}\left( B\right) &=&k\left( B\right) ^{\dagger }k_{\bullet }\left( B\right)
=\gamma ^{\bullet }\left( B^{\dagger }B\right) ^{\dagger }-\gamma ^{\bullet
}\left( B\right) ^{\dagger }B,
\end{eqnarray*}
the property (\ref{3.8}) is fulfilled, $L_{\circ }^{\bullet }j\left(
B\right) L_{\bullet }^{\circ }=\gamma _{\bullet }^{\bullet }\left( B\right) $
with $L_{\circ }^{\bullet }=k_{\bullet }^{*}\left( I\right) =L_{\bullet
}^{\circ \dagger }$, and 
\begin{equation*}
\gamma _{\bullet }\left( B^{\dagger }\right) =k\left( B\right) ^{\dagger
}L_{\bullet }^{\circ }+B^{\dagger }L_{\bullet }^{-}=\gamma ^{\bullet }\left(
B\right) ^{\dagger },
\end{equation*}
where $L_{\bullet }^{-}=\gamma _{\bullet }\left( I\right) ,L_{+}^{\bullet
}=\gamma ^{\bullet }\left( I\right) =L_{\bullet }^{-\dagger }$.

The proof of the implication (ii)$\Rightarrow $(iii) can be also obtained as
in \cite{18} by the explicit construction of $\mathcal{E}$ as $\mathcal{H}%
\oplus \mathcal{H}^{\circ }\oplus \mathcal{H}$ with the indefinite metric
tensor $\mathbf{G}=\left[ G_{\mu \nu }\right] $ given above for $\mu ,\nu
=-,\circ ,+$, and $D=\gamma \left( I\right) $. The unital $\flat $%
-representation $\jmath =\left[ \jmath _\nu ^\mu \right] _{\nu =-,\circ
,+}^{\mu =-,\circ ,+}$ of $\mathcal{B}$ on $\mathcal{E}$ : 
\begin{equation*}
\jmath \left( X^{\dagger }Z\right) =\jmath \left( X\right) ^{\flat }\jmath
\left( Z\right) ,\quad \jmath \left( I\right) =I
\end{equation*}
with $\jmath \left( B\right) ^{\flat }=\mathbf{G}^{-1}\jmath \left( B\right)
^{\dagger }\mathbf{G}=\jmath \left( B^{\dagger }\right) $ is given by the
components 
\begin{equation}
\jmath _{\circ }^{\circ }=j,\quad \jmath _{+}^{\circ }=k,\quad \jmath
_{\circ }^{-}=k^{*},\quad \jmath _{+}^{-}=l,\quad \jmath _{-}^{-}=i=\jmath
_{+}^{+}  \label{3.9}
\end{equation}
and all other $\jmath _\nu ^\mu =0$. The linear operator $\mathbf{L}:%
\mathcal{H}\oplus \mathcal{H}^{\bullet }\rightarrow \mathcal{E}$, where $%
\mathcal{H}^{\bullet }=\mathcal{H}\otimes \mathbb{C}^d$, can be defined by
the components $\left( L^\mu ,L_{\bullet }^\mu \right) $, 
\begin{equation*}
L^{-}=0,\quad L^{\circ }=0,\quad L^{+}=I,\quad L_{\bullet }^{-}=\left(
L_n^{-}\right) ,\quad L_{\bullet }^{\circ }=\left( L_n^{\circ }\right)
,\quad L_{\bullet }^{+}=0,
\end{equation*}
and $\mathbf{L}^{\flat }=\left( 
\begin{array}{ccc}
I & 0 & D \\ 
0 & L_{\circ }^{\bullet } & L_{+}^{\bullet }%
\end{array}
\right) =\mathbf{L}^{\dagger }\mathbf{G}$, where $L_{\circ }^{\bullet
}=L_{\bullet }^{\circ \dagger },L_{+}^{\bullet }=L_{\bullet }^{-\dagger }$.
Then $\mathbf{L}^{\flat }\jmath \mathbf{L=}$ 
\begin{equation*}
\left( 
\begin{array}{cc}
0 & L_{\bullet }^{-} \\ 
0 & L_{\bullet }^{\circ } \\ 
1 & 0%
\end{array}
\right) ^{\flat }\left[ 
\begin{array}{ccc}
i & k^{*} & l \\ 
0 & j & k \\ 
0 & 0 & i%
\end{array}
\right] \left( 
\begin{array}{cc}
0 & L_{\bullet }^{-} \\ 
0 & L_{\bullet }^{\circ } \\ 
1 & 0%
\end{array}
\right) =\left( 
\begin{array}{cc}
l+Di & k^{*}L_{\bullet }^{\circ }+iL_{\bullet }^{-} \\ 
L_{\circ }^{\bullet }k+L_{+}^{\bullet }i & L_{\circ }^{\bullet }jL_{\bullet
}^{\circ }%
\end{array}
\right) =\boldsymbol{\gamma}
\end{equation*}

In order to prove the implication (iii)$\Rightarrow $(iv), it is sufficient
to show that the vectors $\xi =\sum_B\jmath \left( B\right) \mathbf{L}%
\boldsymbol{\eta }_B$ are positive, $\left( \xi |\xi \right) \geq 0$ if $%
\sum_B\boldsymbol{\iota }\left( B\right) \boldsymbol{\eta }_B=\sum_BB\eta
_B=0$. But this follows immediately from the observation $\xi
^{+}=\sum_B\jmath \left( B\right) L^{+}\eta _B=\sum_BB\eta _B=0$ such that
the indefinite metrics (\ref{3.4}) is positive, $\left( \xi |\xi \right)
=\left\| \xi ^{\circ }\right\| ^2\geq 0$ in this case.

The final implication (iv)$\Rightarrow $(i) is obtained as the case $\eta
_I=-\sum_{B\neq I}B\eta _B$ of $\sum_BB\eta _B=0$. 
%TCIMACRO{\TeXButton{End Proof}{\endproof}}%
%BeginExpansion
\endproof%
%EndExpansion

\section{The Structure of the Bounded Filtering Generators.}

The structure (\ref{3.3}) of the form-generator for CP cocycles over $%
\mathcal{B}=\mathcal{B}\left( \mathcal{H}\right) $ is a consequence of the
well known fact that the derivations $k,k^{*}$ of the algebra $\mathcal{B}%
\left( \mathcal{H}\right) $ of all bounded operators on a Hilbert space $%
\mathcal{H}$ are spatial, $k\left( B\right) =j\left( B\right) L-LB,\quad
k^{*}\left( B\right) =L^{\dagger }j\left( B\right) -BL^{\dagger }$, and so 
\begin{equation}
l\left( B\right) =\frac 12\left( L^{\dagger }k\left( B\right) +k^{*}\left(
B\right) L+\left[ B,D\right] \right) +i\left[ H,B\right] ,  \label{4.1}
\end{equation}
where $H^{\dagger }=H$ is a Hermitian operator in $\mathcal{H}$. The
germ-map $\boldsymbol{\gamma }$ whose components are composed (as in (\ref%
{3.3})) into the sums of the components $\varphi _\nu ^\mu \quad $of a CP
matrix map $\boldsymbol{\varphi }:\mathcal{B}\rightarrow \mathcal{B}\otimes 
\mathcal{M}\left( \mathbb{C}^{d+1}\right) $ and left and right
multiplications, are obviously conditionally completely positive with
respect to the representation $\boldsymbol{\iota }$ in (4). As follows from
the dilation theorem in this case, there exists a family $%
L_{-}=L=L_{+},\quad L_n=L_n^{\circ },\quad n=1,...,d$ of linear operators $%
L_\nu :\mathcal{H}\rightarrow \mathcal{H}^{\circ }$, having adjoints $L_\mu
^{\dagger }:\mathcal{H}^{\circ }\rightarrow \mathcal{H}$ such that $\varphi
_\nu ^\mu \left( B\right) =L_\mu ^{\dagger }j\left( B\right) L_\nu $.

The next theorem proves that these structural conditions which are
sufficient for complete positivity of the cocycles, given by the equation (%
\ref{2.2}), are also necessary if the germ-map $\boldsymbol{\gamma }$ is
w*-continuous on an operator algebra $\mathcal{B}$. Thus the equation (\ref%
{2.2}) for a completely positive quantum cocycle with bounded stochastic
derivatives has the following general form 
\begin{equation*}
\mathrm{d}\phi _t\left( B\right) +\phi _t\left( K^{\dagger }B+BK-L^{\dagger
}j\left( B\right) L\right) \mathrm{d}t=\sum_{m,n=1}^d\phi _t\left(
L_m^{\dagger }j\left( B\right) L_n-B\delta _n^m\right) \mathrm{d}\Lambda _m^n
\end{equation*}
\begin{equation}
+\sum_{m=1}^d\phi _t\left( L_m^{\dagger }j\left( B\right) L-K_m^{\dagger
}B\right) \mathrm{d}\Lambda _m^{+}+\sum_{n=1}^d\phi _t\left( L^{\dagger
}j\left( B\right) L_n-BK_n\right) \mathrm{d}\Lambda _{-}^n,  \label{4.2}
\end{equation}
generalising the Lindblad form \cite{17}, for the norm-continuous semigroups
of completely positive maps. The quantum stochastic submartingale $M_t=\phi
_t\left( I\right) $ is defined by the integral 
\begin{equation*}
M_t+\int_0^t\phi _s\left( D\right) \mathrm{d}s=I+\int_0^t\sum_{m,n}^d\phi
_s\left( L_n^{\dagger }L_m-\delta _m^n\right) \mathrm{d}\Lambda _n^m
\end{equation*}
\begin{equation}
+\int_0^t\sum_{m=1}^d\phi _s\left( L_m^{\dagger }L-K_m^{\dagger }\right) 
\mathrm{d}\Lambda _m^{+}+\int_0^t\sum_{n=1}^d\phi _s\left( L^{\dagger
}L_n-K_n\right) \mathrm{d}\Lambda _{-}^n.  \label{4.3}
\end{equation}
If the space $\mathcal{K}$ can be embedded into the direct sum $\mathcal{H}%
\otimes \mathbb{C}^d=\mathcal{H}\oplus ...\oplus \mathcal{H}$ of $d$ copies
of the initial Hilbert space $\mathcal{H}$ such that $j\left( B\right)
=\left( B\delta _n^m\right) $, this equation can be resolved in the form $%
\phi _t\left( B\right) =F_t^{\dagger }BF_t$, where $F=\left( F_t\right)
_{t>0}$ is an (unbounded) cocycle in the tensor product $\mathcal{H}\otimes 
\mathcal{F}$ with Fock space $\mathcal{F}$ over the Hilbert space $\mathbb{C}%
^d\otimes L^2\left( \mathbb{R}_{+}\right) $ of the quantum noise of
dimensionality $d$. The cocycle $F$ satisfies the quantum stochastic
equation 
\begin{equation}
\mathrm{d}F_t+KF_t\mathrm{d}t=\sum_{i,n=1}^d\left( L_n^i-I\delta _n^i\right)
F_t\mathrm{d}\Lambda _i^n+\sum_{i=1}^dL^iF_t\mathrm{d}\Lambda
_i^{+}-\sum_{n=1}^dK_nF_t\mathrm{d}\Lambda _{-}^n,\qquad  \label{4.4}
\end{equation}
where $L_n^{i\text{ }}$ and $L^i$ are the operators in $\mathcal{H}$,
defining 
\begin{eqnarray}
\varphi _n^m\left( B\right) &=&\sum_{i=1}^dL_m^{i\dagger }BL_n^i,\qquad
\varphi \left( B\right) =\sum_{i=1}^dL^{i\dagger }BL^i  \label{4.5} \\
\varphi ^m\left( B\right) &=&\sum_{i=1}^dL_m^{i\dagger }BL^i,\qquad \varphi
_n\left( B\right) =\sum_{i=1}^dL^{i\dagger }BL_n^i\qquad  \notag
\end{eqnarray}
with $\sum_{i=1}^dL^{i\dagger }L^i=K+K^{\dagger }$ if $M_t$ is a martingale (%
$\leq K+K^{\dagger }$if submartingale) .\allowbreak

\begin{theorem}
Let the germ-maps $\boldsymbol{\gamma }$ of the quantum stochastic cocycle $%
\phi $ over a von-Neumann algebra $\mathcal{B}$ be w*-continuous and
bounded: 
\begin{equation}
\left\| \gamma \right\| <\infty ,\qquad \left\| \gamma _{\bullet }\right\|
=\left( \sum_{n=1}^d\left\| \gamma _n\right\| ^2\right) ^{\frac 12}=\left\|
\gamma ^{\bullet }\right\| <\infty ,\qquad \left\| \gamma _{\bullet
}^{\bullet }\right\| =\left\| \gamma _{\bullet }^{\bullet }\left( I\right)
\right\| <\infty ,  \label{4.6}
\end{equation}
where $\left\| \gamma \right\| =\sup \left\{ \left\| \gamma \left( B\right)
\right\| :\left\| B\right\| <1\right\} ,\quad \left\| \gamma _{\bullet
}^{\bullet }\left( I\right) \right\| =\sup \left\{ \left\langle \eta
^{\bullet },\gamma _{\bullet }^{\bullet }\left( I\right) \eta ^{\bullet
}\right\rangle \left| \left\| \eta ^{\bullet }\right\| <1\right. \right\} $
and $\phi _t$ be a CP cocycle, satisfying equation (\ref{2.2}) with $\phi
_0\left( B\right) =B$ and normalized to a submartingale (martingale). Then
they have the form (\ref{3.3}) written as 
\begin{equation}
\boldsymbol{\gamma }\left( B\right) =\boldsymbol{\varphi }\left( B\right) -%
\boldsymbol{\iota }\left( B\right) \boldsymbol{K}-\boldsymbol{K}^{\dagger }%
\boldsymbol{\iota }\left( B\right)  \label{4.7}
\end{equation}
with $\varphi =\varphi _{+}^{-},\quad \varphi ^m=\varphi _{+}^m,\quad
\varphi _n=\varphi _n^{-}$ and $\varphi _n^m=\gamma _n^m$, composing a
bounded CP map. 
\begin{equation}
\boldsymbol{\varphi }=\left( 
\begin{array}{cc}
\varphi & \varphi _{\bullet } \\ 
\varphi ^{\bullet } & \varphi _{\bullet }^{\bullet }%
\end{array}
\right) ,\quad and\quad \boldsymbol{K}=\left( 
\begin{array}{cc}
K & K_{\bullet } \\ 
K^{\bullet } & K_{\bullet }^{\bullet }%
\end{array}
\right)  \label{4.8}
\end{equation}
with arbitrary $K^{\bullet },K_{\bullet }^{\bullet }$, and $K+K^{\dagger
}\geq \varphi \left( I\right) $. The equation (\ref{4.2}) has the unique CP
solution , satisfying the condition $\phi _s\left( I\right) \leq \epsilon _s%
\left[ \phi _t\left( I\right) \right] $ for all $s<t$ ($\phi _s\left(
I\right) =\epsilon _s\left[ \phi _t\left( I\right) \right] $ if $%
K+K^{\dagger }=\varphi \left( I\right) $).
\end{theorem}

%TCIMACRO{\TeXButton{Proof}{\proof}}%
%BeginExpansion
\proof%
%EndExpansion
The structure (\ref{4.7}) for the CP component $\gamma _{\bullet }^{\bullet
} $ was obtained as a part of the dilation theorem in the Stinespring form $%
\gamma _{\bullet }^{\bullet }\left( B\right) =L_{\bullet }^{\dagger }j\left(
B\right) L_{\bullet }=\varphi _{\bullet }^{\bullet }\left( B\right) $, where 
$L_{\bullet }=L_{\bullet }^{\circ }$. In order to obtain the structure (\ref%
{4.7}) for the bounded germ-maps $\gamma _{\bullet }$ and $\gamma
_{}^{\bullet }$, we can take into account the spatial structure $k\left(
B\right) =j\left( B\right) L-LB$ of a bounded $\left( j,i\right) $%
-derivation for a von-Neumann algebra $\mathcal{B}$ with respect to a normal
representation $j$ of $\mathcal{B}$ and $i\left( B\right) =B$. Then 
\begin{equation*}
\gamma _{\bullet }\left( B\right) =k^{*}\left( B\right) L_{\bullet }^{\circ
}+BL_{\bullet }^{-}=L^{\dagger }j\left( B\right) L_{\bullet }^{\circ
}-B\left( L^{\dagger }L_{\bullet }^{\circ }-L_{\bullet }^{-}\right)
=L_{+}^{\dagger }j\left( B\right) L_{\bullet }-BK_{\bullet },
\end{equation*}
where $L_{+}=L,K_{\bullet }=L^{\dagger }L_{\bullet }^{\circ }-L_{\bullet
}^{-}$. Hence $\gamma _{\bullet }\left( B\right) =\varphi _{\bullet }\left(
B\right) -BK_{\bullet },\gamma ^{\bullet }\left( B\right) =\varphi ^{\bullet
}\left( B\right) -K^{\bullet }B=\gamma _{\bullet }^{*}\left( B\right) $,
where $K^{\bullet }=K_{\bullet }^{\dagger },\varphi ^{\bullet }\left(
B\right) =L_{\bullet }^{\dagger }j\left( B\right) L=\varphi _{\bullet
}^{*}\left( B\right) $, such that the matrix-map $\boldsymbol{\varphi }%
\left( B\right) =\left( L^\mu j\left( B\right) L_\nu \right) _{\nu
=+,\bullet }^{\mu =-,\bullet }$ with $L^{-}=L^{\dagger },L^{\bullet
}=L_{\bullet }^{\dagger }$ is CP. Taking into account the form (\ref{4.1})
of the coboundary $l\left( B\right) =\gamma \left( B\right) -DB$ which is
due to the spatial form $\left[ iH+\frac 12D,B\right] $ of the bounded
derivation $l\left( B\right) -\frac 12\left( L^{\dagger }k\left( B\right)
+k^{*}\left( B\right) L\right) $ on $\mathcal{B}$, one can obtain the
representation 
\begin{equation*}
\gamma \left( B\right) =\frac 12\left( L^{\dagger }k\left( B\right)
+k^{*}\left( B\right) L+DB+BD\right) +i\left[ H,B\right] =\varphi \left(
B\right) -BK-K^{\dagger }B,
\end{equation*}
where $\varphi \left( B\right) =L^{\dagger }j\left( B\right) L$, $K=iH+\frac
12\left( L^{\dagger }L-D\right) $.

The existence and uniqueness of the solutions $\phi _t\left( B\right) $ to
the quantum stochastic equations (\ref{2.2}) with the bounded generators $%
\lambda _\nu ^\mu \left( B\right) =\gamma _\nu ^\mu \left( B\right) -B\delta
_\nu ^\mu $ and the initial conditions $\phi _0\left( B\right) =B$ in an
operator algebra $\mathcal{B}$ was proved in \cite{20}. The positivity of
the solutions in the case of the equation (\ref{4.2}), corresponding to the
conditionally positive germ-function (\ref{4.7}), can be obtained by the
iteration 
\begin{equation*}
\phi _t^{\left( n+1\right) }\left( B\right) =V_t^{\dagger }BV_t+\int_0^t\phi
_s^{\left( n\right) }\left( \beta _\nu ^\mu \left( V_t^{\dagger }\left(
s\right) BV_t\left( s\right) \right) \right) \mathrm{d}\Lambda _\mu ^\nu
,\quad \phi _t^{\left( 0\right) }\left( B\right) =B
\end{equation*}
of the quantum stochastic integral equation 
\begin{equation}
\phi _t\left( B\right) =V_t^{\dagger }BV_t+\int_0^t\phi _s\left( \beta _\nu
^\mu \left( V_t^{\dagger }\left( s\right) BV_t\left( s\right) \right)
\right) \mathrm{d}\Lambda _\mu ^\nu ,  \label{4.9}
\end{equation}
with $\beta _\nu ^\mu \left( B\right) =\varphi _\nu ^\mu \left( B\right)
-B\delta _\nu ^\mu $. Here $V_t=V_t\left( s\right) V_s$ with $V_t\left(
s\right) =T_s^{\dagger }V_{t-s}T_s$ shifted by the co-isometry $T_s$ in $%
\mathcal{D}$, is the vector cocycle, resolving the quantum stochastic
differential equation 
\begin{equation}
\mathrm{d}V_t+KV_t\mathrm{d}t+\sum_{m=1}^dK_mV_t\mathrm{d}\Lambda _{-}^m=0
\label{4.10}
\end{equation}
with the initial condition $V_0=I$ in $\mathcal{H}$. The equivalence of (\ref%
{4.2}) and (\ref{4.9}), (\ref{4.10}) is verified by direct differentiation
of (\ref{4.9}). In order to prove the complete positivity of this solution,
one should write down the corresponding iteration 
\begin{equation*}
\phi _t^{\left( n+1\right) }\left( f^{\bullet },B,h^{\bullet }\right)
=V_t^{\dagger }BV_t+\int_0^t\boldsymbol{f}\left( s\right) ^{\dagger }\phi
_s^{\left( n\right) }\left( f^{\bullet },\boldsymbol{\varphi }\left(
V_t^{\dagger }\left( s\right) BV_t\left( s\right) \right) ,h^{\bullet
}\right) \boldsymbol{h}\left( s\right) \mathrm{d}s,
\end{equation*}
of the ordinary integral equation for the operator-valued kernels of
coherent vectors, defined in (\ref{2.7}). Here $\boldsymbol{g}\left(
s\right) =1\oplus g^{\bullet }\left( s\right) $ such that 
\begin{equation*}
\sum_{X,Z}\sum_{f,h}\left\langle \xi _X^f\right| \left. \phi _t\left(
f^{\bullet },X^{\dagger }Z,h^{\bullet }\right) \xi _Z^h\right\rangle
=\sum_{X,Z}\left\langle XV_t\eta _X|ZV_t\eta _Z\right\rangle
\end{equation*}
\begin{equation*}
+\int_0^t\sum_{X,Z}\sum_{f,h}\left\langle \boldsymbol{\eta }_X^f\left(
s\right) |\phi _s\left( f^{\bullet },\boldsymbol{\varphi }\left( X^{\dagger
}Z\right) ,h^{\bullet }\right) \boldsymbol{\eta }_Z^h\left( s\right)
\right\rangle ,
\end{equation*}
where $\eta _B=\sum_g\xi _B^g,\quad \boldsymbol{\eta }_B^g\left( s\right)
=\sum_{g\left( s\right) }\xi _B^g\otimes \boldsymbol{g}\left( s\right) $.
Then the CP property for $\phi _t^{\left( n\right) }$, immediately follows
from the CP property of $\phi _s^{\left( n-1\right) },s<t$ and of $%
\boldsymbol{\varphi }$. The direct iteration of this integral recursion with
the initial CP condition $\phi _t^{\left( 0\right) }\left( B\right) =B$
gives at the limit $n\rightarrow \infty $ the minimal CP solution in the
form of sum of n-tupol CP integrals on the interval $\left[ 0,t\right] $. 
%TCIMACRO{\TeXButton{End Proof}{\endproof}}%
%BeginExpansion
\endproof%
%EndExpansion

\end{document}